\documentclass[12pt,fleqn]{article}

\usepackage{amscd,amsmath,amssymb,amsthm}
\usepackage{geometry} 
\usepackage[pdfpagemode=UseNone,pdfstartview=FitH,hypertex]{hyperref}

\newcommand{\Ap}[1][]{A_p\, #1}

\newcommand{\Bp}[1][]{B_p\, #1}
\newcommand{\bdry}[1]{\partial #1}
\newcommand{\A}{{\cal A}}
\newcommand{\D}{{\cal D}}
\newcommand{\F}{{\cal F}}
\newcommand{\closure}[1]{\overline{#1}}

\newcommand{\dint}{\ds{\int}}
\newcommand{\dist}[2]{\text{dist}\, (#1,#2)}
\newcommand{\dnorm}[2][]{\left\|#2\right\|_{#1}^\ast}
\newcommand{\ds}[1]{\displaystyle #1}
\newcommand{\dualp}[3][]{\left(#2,#3\right)_{#1}}
\newcommand{\eps}{\varepsilon}

\newcommand{\hquad}{\hspace{0.08in}}
\newcommand{\incl}{\hookrightarrow}

\newcommand{\loc}{\text{loc}}
\newcommand{\M}{{\cal M}}
\newcommand{\N}{\mathbb N}
\newcommand{\norm}[2][]{\left\|#2\right\|_{#1}}

\renewcommand{\o}{\text{o}}
\newcommand{\PS}[1]{$(\text{PS})_{#1}$}
\newcommand{\pnorm}[2][]{\if #1'' \left|#2\right|_p \else \left|#2\right|_{#1} \fi}

\newcommand{\R}{\mathbb R}
\newcommand{\RP}{\R \text{P}}
\newcommand{\restr}[2]{\left.#1\right|_{#2}}
\newcommand{\seq}[1]{\left(#1\right)}
\newcommand{\set}[1]{\left\{#1\right\}}
\newcommand{\wto}{\rightharpoonup}
\newcommand{\Z}{\mathbb Z}

\DeclareMathOperator{\divg}{div}
\DeclareMathOperator{\sgn}{sgn}

\newenvironment{enumroman}{\begin{enumerate}
\renewcommand{\theenumi}{$(\roman{enumi})$}
\renewcommand{\labelenumi}{$(\roman{enumi})$}}{\end{enumerate}}

\newenvironment{properties}[1]{\begin{enumerate}
\renewcommand{\theenumi}{$(#1_\arabic{enumi})$}
\renewcommand{\labelenumi}{$(#1_\arabic{enumi})$}}{\end{enumerate}}

\newtheorem{corollary}{Corollary}[section]
\newtheorem{lemma}[corollary]{Lemma}
\newtheorem{proposition}[corollary]{Proposition}
\newtheorem{theorem}[corollary]{Theorem}

\theoremstyle{definition}
\newtheorem{definition}[corollary]{Definition}

\theoremstyle{remark}
\newtheorem{example}[corollary]{Example}

\numberwithin{equation}{section}

\title{\bf New multiplicity results for critical $p$-Laplacian problems}
\author{\bf Carlo Mercuri\\
Department of Mathematics\\
Swansea University\\
Swansea, SA1 8EN, UK\\
\em C.Mercuri@swansea.ac.uk\\
[\bigskipamount]
\bf Kanishka Perera\\
Department of Mathematical Sciences\\
Florida Institute of Technology\\
Melbourne, FL 32901, USA\\
\em kperera@fit.edu}
\date{}

\begin{document}

\maketitle

\begin{abstract}
We prove new multiplicity results for the Br{\'e}zis-Nirenberg problem for the $p$-Laplacian. Our proofs are based on a new abstract critical point theorem involving the $\Z_2$-cohomological index that requires less compactness than the \PS{} condition.
\end{abstract}

\section{Introduction and statement of results}

Consider the problem
\begin{equation} \label{1}
\left\{\begin{aligned}
- \Delta_p\, u & = \lambda\, |u|^{p-2}\, u + |u|^{p^\ast - 2}\, u && \text{in } \Omega\\[10pt]
u & = 0 && \text{on } \bdry{\Omega},
\end{aligned}\right.
\end{equation}
where $\Omega$ is a smooth bounded domain in $\R^N,\, N \ge 2$, $\Delta_p\, u = \divg (|\nabla u|^{p-2}\, \nabla u)$ is the $p$-Laplacian of $u$, $1 < p < N$, $p^\ast = Np/(N - p)$ is the critical Sobolev exponent, and $\lambda > 0$ is a parameter. Existence and multiplicity of nontrivial solutions to this problem has been widely studied beginning with the celebrated paper of Br{\'e}zis and Nirenberg \cite{MR709644} (see, e.g., [1--8, 11--15, 17--24, 27, 28, 30--35, 37, 40, 42, 43, 46, 48]). 
In particular, the following multiplicity results are known in the semilinear case $p = 2$:
\begin{enumroman}
\item If $N \ge 7$, then problem \eqref{1} has infinitely many solutions for all $\lambda > 0$ (see Devillanova and Solimini \cite{MR1919704}).

\item If $4 \le N \le 6$ and $0 < \lambda < \lambda_1$, where $\lambda_1 > 0$ is the first Dirichlet eigenvalue of $- \Delta$ on $\Omega$, then problem \eqref{1} has $(N + 2)/2$ distinct pairs of nontrivial solutions (see Devillanova and Solimini \cite{MR1966256} and Clapp and Weth \cite{MR2122698}).

\item If $4 \le N \le 6$ and $\lambda > \lambda_1$ is not an eigenvalue, then problem \eqref{1} has $(N + 1)/2$ distinct pairs of nontrivial solutions (see Clapp and Weth \cite{MR2122698}).

\item If $N = 5 \text{ or } 6$ and $\lambda \ge \lambda_1$, then problem \eqref{1} has $(N + 1)/2$ distinct pairs of nontrivial solutions (see Chen et al.\! \cite{MR2926297}).

\item If $N = 4$ and $\lambda \ge \lambda_1$ is an eigenvalue of multiplicity $m \ge 1$, then problem \eqref{1} has $(N - m + 1)/2$ distinct pairs of nontrivial solutions (see Clapp and Weth \cite{MR2122698}).
\end{enumroman}
In the general case $1 < p < N$, the following multiplicity results are known:
\begin{enumroman}
\item If $N > p^2 + p$, then problem \eqref{1} has infinitely many solutions for all $\lambda > 0$ (see Cao et al.\! \cite{MR2885967}).

\item If $p^2 \le N \le p^2 + p$ and $0 < \lambda < \lambda_1$, where $\lambda_1 > 0$ is the first Dirichlet eigenvalue of $- \Delta_p$ on $\Omega$, then problem \eqref{1} has infinitely many solutions (see He et al.\! \cite{MR4098039}).
\end{enumroman}
The purpose of the present paper is to prove some multiplicity results for the case where $p^2 \le N \le p^2 + p$ and $\lambda \ge \lambda_1$ similar to those of Clapp and Weth \cite{MR2122698} and Chen et al.\! \cite{MR2926297}. However, the arguments in \cite{MR2122698} and \cite{MR2926297} are based on the relative equivariant Lusternik-Schnirelmann category and the Krasnoselskii genus, respectively, and make essential use of the fact that the Laplacian is a linear operator and therefore has linear eigenspaces. These arguments do not extend to the $p$-Laplacian, which is a nonlinear operator and hence lacks linear eigenspaces. Our proofs will make use of a new abstract critical point theorem based on the $\Z_2$-cohomological index that we will prove in the next section. This theorem has no hypotheses involving linear subspaces and is of independent interest.

To state our multiplicity results, recall that solutions of problem \eqref{1} coincide with critical points of the energy functional
\[
E(u) = \frac{1}{p} \int_\Omega |\nabla u|^p\, dx - \frac{\lambda}{p} \int_\Omega |u|^p\, dx - \frac{1}{p^\ast} \int_\Omega |u|^{p^\ast} dx, \quad u \in W^{1,\,p}_0(\Omega),
\]
and that eigenvalues of the asymptotic eigenvalue problem
\[
\left\{\begin{aligned}
- \Delta_p\, u & = \lambda\, |u|^{p-2}\, u && \text{in } \Omega\\[10pt]
u & = 0 && \text{on } \bdry{\Omega}
\end{aligned}\right.
\]
coincide with critical values of the functional
\[
\Psi(u) = \frac{1}{\dint_\Omega |u|^p\, dx}, \quad u \in S = \set{u \in W^{1,\,p}_0(\Omega) : \int_\Omega |\nabla u|^p\, dx = 1}.
\]
Denote by $\F$ the class of symmetric subsets of $S$ and by $i(M)$ the $\Z_2$-cohomological index of $M \in \F$ (see Definition \ref{Definition 1}), let $\F_k = \set{M \in \F : i(M) \ge k}$, and set
\[
\lambda_k = \inf_{M \in \F_k}\, \sup_{u \in M}\, \Psi(u), \quad k \ge 1.
\]
Then $\lambda_1 = \inf \Psi(S) > 0$ is the first eigenvalue and $\lambda_1 < \lambda_2 \le \cdots$ is an unbounded sequence of eigenvalues. Let
\begin{equation} \label{20}
S_{N,\,p} = \inf_{u \in \D^{1,\,p}(\R^N) \setminus \set{0}}\, \frac{\dint_{\R^N} |\nabla u|^p\, dx}{\left(\dint_{\R^N} |u|^{p^\ast} dx\right)^{p/p^\ast}}
\end{equation}
denote the best Sobolev constant. Our first multiplicity result is the following.

\begin{theorem} \label{Theorem 5}
Let $N \ge p^2$.
\begin{enumroman}
\item \label{Theorem 5.i} If $0 < \lambda < \lambda_1$ or $\lambda_k < \lambda < \lambda_{k+1}$ for some $k \ge 1$, then problem \eqref{1} has $N/2$ distinct pairs of nontrivial solutions.
\item \label{Theorem 5.ii} If $\lambda = \lambda_1$ and $N \ge 3$, then problem \eqref{1} has $(N - 1)/2$ distinct pairs of nontrivial solutions.
\item \label{Theorem 5.iii} If $\lambda_{k-m} < \lambda = \lambda_{k-m+1} = \cdots = \lambda_k < \lambda_{k+1}$ for some $k > m \ge 1$ and $N \ge m + 2$, then problem \eqref{1} has $(N - m)/2$ distinct pairs of nontrivial solutions.
\end{enumroman}
These solutions satisfy
\begin{equation} \label{23}
0 < E(u) < \frac{2}{N}\, S_{N,\,p}^{N/p}.
\end{equation}
\end{theorem}

Eigenvalues based on the cohomological index were first introduced in Perera \cite{MR1998432} (see also Perera and Szulkin \cite{MR2153141}). A complete description of the spectrum of the $p$-Laplacian is not known when $p \ne 2$, and $\seq{\lambda_k}$ may not be a complete list of eigenvalues. However, Theorem \ref{Theorem 5} \ref{Theorem 5.i} gives $N/2$ distinct pairs of nontrivial solutions of problem \eqref{1} whenever $\lambda > 0$ is not an eigenvalue from the sequence $\seq{\lambda_k}$.

Our second multiplicity result makes no references to the spectrum.

\begin{theorem} \label{Theorem 6}
If $N^2/(N + 1) > p^2$, then problem \eqref{1} has $N/2$ distinct pairs of nontrivial solutions satisfying
\begin{equation} \label{24}
0 < E(u) < \frac{2}{N}\, S_{N,\,p}^{N/p}
\end{equation}
for all $\lambda > 0$.
\end{theorem}

As is usually the case with problems of critical growth, the energy functional $E$ associated with problem \eqref{1} does not satisfy the \PS{c} condition for all $c \in \R$. However, it has certain weaker compactness properties (see Theorem \ref{PScondition}). We will first prove an abstract multiplicity result that only assumes these weaker conditions, and apply it to prove Theorems \ref{Theorem 5} and \ref{Theorem 6}.

To state our abstract results, let $D$ be a Banach space and let $W$ be closed linear subspace of $D$. For $A \subset D$ and $\delta > 0$, we set
\[
N_\delta(A) = \set{u \in W : \dist{u}{A} \le \delta}.
\]
Let $E$ be an even $C^1$-functional on $W$ such that $E(0) = 0$. For $c \in \R$, let
\[
K_c = \set{u \in W : E'(u) = 0,\, E(u) = c}
\]
be the set of critical points of $E$ at the level $c$. We assume that $E$ has the following compactness properties:
\begin{properties}{C}
\item \label{C1} there exists $c^\ast > 0$ such that $E$ satisfies the \PS{c} condition for all $c \in (0,c^\ast)$,
\item \label{C2} there exist $b > c^\ast$ and for each $c \in [c^\ast,b)$ a set $M_c \subset D \setminus \set{0}$ such that
    \begin{enumerate}
    \item[({\em i})] every \PS{c} sequence $\seq{u_n}$ has either a subsequence that converges strongly to a point in $K_c$, or a renamed subsequence that converges weakly to a point in $K_{c-c^\ast}$ and satisfies
        \[
        \dist{u_n}{M_c} \to 0 \quad \text{or} \quad \dist{u_n}{-M_c} \to 0,
        \]
    \item[({\em ii})] $N_\delta(M_c) \cap N_\delta(-M_c) = \emptyset$ for all sufficiently small $\delta > 0$.
    \end{enumerate}
\end{properties}
Let $S = \set{u \in W : \norm{u} = 1}$ be the unit sphere in $W$, let $\pi : W \setminus \set{0} \to S,\, u \mapsto u/\norm{u}$ be the radial projection onto $S$, and let $S^N = \set{x \in \R^{N+1} : |x| = 1}$ be the unit sphere in $\R^{N+1}$. We will prove the following theorem in the next section.

\begin{theorem} \label{Theorem 3}
Assume \ref{C1} and \ref{C2}. Let $B_0$ and $C_0$ be symmetric subsets of $S$ such that $C_0$ is compact, $B_0$ is closed, and
\begin{equation} \label{11}
i(C_0) \ge k - m, \qquad i(S \setminus B_0) \le k
\end{equation}
for some $k \ge m \ge 0$. Assume that there exist an odd continuous map $\varphi : S^N \to S \setminus C_0,\, N \ge m + 2$ and $R > r > 0$ such that, setting
\begin{gather*}
A_0 = \begin{cases}
\varphi(S^N) & \text{if } C_0 = \emptyset\\[7.5pt]
\set{\pi((1 - t)\, v + tw) : v \in C_0,\, w \in \varphi(S^N),\, t \in [0,1]} & \text{if } C_0 \ne \emptyset,
\end{cases}\\[10pt]
A = \set{Ru : u \in A_0}, \qquad B = \set{ru : u \in B_0}, \qquad X = \set{tu : u \in A,\, t \in [0,1]},
\end{gather*}
we have
\[
\sup_A\, E \le 0 < \inf_B\, E, \qquad \sup_X\, E < b.
\]
Then $E$ has $(N - m)/2$ distinct pairs of critical points at levels in $(0,b)$.
\end{theorem}

We will apply this theorem to a class of abstract $p$-Laplacian equations that includes problem \eqref{1} as a special case. Assume that $W$ is uniformly convex and let $(W^\ast,\dnorm{\,\cdot\,})$ be its dual with duality pairing $\dualp{\cdot}{\cdot}$. Recall that $f \in C(W,W^\ast)$ is a potential operator if there is a functional $F \in C^1(W,\R)$, called a potential for $f$, such that $F' = f$. Consider the nonlinear operator equation
\begin{equation} \label{12}
\Ap[u] = \lambda \Bp[u] + f(u)
\end{equation}
in $W^\ast$, where $\Ap, \Bp, f \in C(W,W^\ast)$ are potential operators satisfying the following assumptions, and $\lambda \in \R$ is a parameter:
\begin{enumerate}
\item[$(A_1)$] $\Ap$ is $(p - 1)$-homogeneous and odd for some $p \in (1,\infty)$: $\Ap[(tu)] = |t|^{p-2}\, t\, \Ap[u]$ for all $u \in W$ and $t \in \R$,
\item[$(A_2)$] $\dualp{\Ap[u]}{v} \le \norm{u}^{p-1} \norm{v}$ for all $u, v \in W$, and equality holds if and only if $\alpha u = \beta v$ for some constants $\alpha, \beta \ge 0$, not both zero (in particular, $\dualp{\Ap[u]}{u} = \norm{u}^p$ for all $u \in W$),
\item[$(B_1)$] $\Bp$ is $(p - 1)$-homogeneous and odd: $\Bp[(tu)] = |t|^{p-2}\, t\, \Bp[u]$ for all $u \in W$ and $t \in \R$,
\item[$(B_2)$] $\dualp{\Bp[u]}{u} > 0$ for all $u \in W \setminus \set{0}$, and $\dualp{\Bp[u]}{v} \le \dualp{\Bp[u]}{u}^{(p-1)/p} \dualp{\Bp[v]}{v}^{1/p}$ for all $u, v \in W$,
\item[$(B_3)$] $\Bp$ is a compact operator,
\item[$(F_1)$] the potential $F$ of $f$ with $F(0) = 0$ satisfies $F(u) = \o(\norm{u}^p)$ as $u \to 0$,
\item[$(F_2)$] $\ds{\lim_{t \to + \infty} \frac{F(tu)}{t^p}} = + \infty$ uniformly on compact subsets of $W \setminus \set{0}$.
\end{enumerate}
Solutions of equation \eqref{12} coincide with critical points of the $C^1$-functional
\[
E(u) = I_p(u) - \lambda J_p(u) - F(u), \quad u \in W,
\]
where
\[
I_p(u) = \frac{1}{p} \dualp{\Ap[u]}{u} = \frac{1}{p} \norm{u}^p, \qquad J_p(u) = \frac{1}{p} \dualp{\Bp[u]}{u}
\]
are the potentials of $\Ap$ and $\Bp$ satisfying $I_p(0) = 0 = J_p(0)$, respectively (see Perera et al.\! \cite[Proposition 1.2]{MR2640827}). Moreover, eigenvalues of the asymptotic eigenvalue problem
\[
\Ap[u] = \lambda \Bp[u]
\]
coincide with critical values of the $C^1$-functional
\[
\Psi(u) = \frac{1}{\dualp{\Bp[u]}{u}}, \quad u \in S,
\]
where $S = \set{u \in W : \dualp{\Ap[u]}{u} = 1}$ is the unit sphere in $W$. Denote by $\F$ the class of symmetric subsets of $S$, let $\F_k = \set{M \in \F : i(M) \ge k}$, and set
\[
\lambda_k = \inf_{M \in \F_k}\, \sup_{u \in M}\, \Psi(u), \quad k \ge 1.
\]
Then $\lambda_1 = \inf \Psi(S) > 0$ is the first eigenvalue and $\lambda_1 \le \lambda_2 \le \cdots$ is an unbounded sequence of eigenvalues. Moreover, denoting by $\Psi^a = \set{u \in S : \Psi(u) \le a}$ (resp. $\Psi_a = \set{u \in S : \Psi(u) \ge a}$) the sublevel (resp. superlevel) sets of $\Psi$, we have
\begin{equation} \label{13}
i(S \setminus \Psi_{\lambda_{k+1}}) \le k \le i(\Psi^{\lambda_k})
\end{equation}
(see Perera et al.\! \cite[Theorem 4.6]{MR2640827}). We have the following multiplicity result for the equation \eqref{12}.

\begin{theorem} \label{Theorem 4}
Suppose $(A_1)$--$(F_2)$ hold and $E$ satisfies \ref{C1} and \ref{C2}. Assume that $\lambda < \lambda_{k+1}$, there exists a compact symmetric subset $C_0$ of $S$ with $i(C_0) \ge k - m$ for some $0 \le m \le k$, and there exists an odd continuous map $\varphi : S^N \to S \setminus C_0,\, N \ge m + 2$ such that
\begin{equation} \label{18}
\sup_{w \in \varphi(S^N),\, t \ge 0}\, E(tw) < b \quad \text{if } C_0 = \emptyset
\end{equation}
and
\begin{equation} \label{19}
\sup_{v \in C_0,\, w \in \varphi(S^N),\, s, t \ge 0}\, E(sv + tw) < b \quad \text{if } C_0 \ne \emptyset.
\end{equation}
Then equation \eqref{12} has $(N - m)/2$ distinct pairs of nontrivial solutions satisfying
\[
0 < E(u) < b.
\]
\end{theorem}

We will prove this theorem in the next section, and apply it to prove Theorems \ref{Theorem 5} and \ref{Theorem 6} in Section \ref{Proofs}.

\section{Abstract multiplicity results}

In this section we prove Theorems \ref{Theorem 3} and \ref{Theorem 4}. We begin by recalling the definition and some properties of the $\Z_2$-cohomological index of Fadell and Rabinowitz \cite{MR0478189}.

\begin{definition} \label{Definition 1}
Let $W$ be a Banach space and let $\A$ denote the class of symmetric subsets of $W \setminus \set{0}$. For $A \in \A$, let $\overline{A} = A/\Z_2$ be the quotient space of $A$ with each $u$ and $-u$ identified, let $f : \overline{A} \to \RP^\infty$ be the classifying map of $\overline{A}$, and let $f^\ast : H^\ast(\RP^\infty) \to H^\ast(\overline{A})$ be the induced homomorphism of the Alexander-Spanier cohomology rings. The cohomological index of $A$ is defined by
\[
i(A) = \begin{cases}
0 & \text{if } A = \emptyset\\[7.5pt]
\sup \set{m \ge 1 : f^\ast(\omega^{m-1}) \ne 0} & \text{if } A \ne \emptyset,
\end{cases}
\]
where $\omega \in H^1(\RP^\infty)$ is the generator of the polynomial ring $H^\ast(\RP^\infty) = \Z_2[\omega]$.
\end{definition}

\begin{example}
The classifying map of the unit sphere $S^N$ in $\R^{N+1},\, N \ge 0$ is the inclusion $\RP^N \incl \RP^\infty$, which induces isomorphisms on the cohomology groups $H^q$ for $q \le N$, so $i(S^N) = N + 1$.
\end{example}

The following proposition summarizes the basic properties of this index.

\begin{proposition}[\cite{MR0478189}] \label{Proposition 2}
The index $i : \A \to \N \cup \set{0,\infty}$ has the following properties:
\begin{properties}{i}
\item \label{Proposition 2.i} Definiteness: $i(A) = 0$ if and only if $A = \emptyset$.
\item \label{Proposition 2.ii} Monotonicity: If there is an odd continuous map from $A$ to $B$ (in particular, if $A \subset B$), then $i(A) \le i(B)$. Thus, equality holds when the map is an odd homeomorphism.
\item Dimension: $i(A) \le \dim W$.
\item Continuity: If $A$ is closed, then there is a closed neighborhood $N \in \A$ of $A$ such that $i(N) = i(A)$. When $A$ is compact, $N$ may be chosen to be a $\delta$-neighborhood $N_\delta(A) = \set{u \in W : \dist{u}{A} \le \delta}$.
\item Subadditivity: If $A$ and $B$ are closed, then $i(A \cup B) \le i(A) + i(B)$.
\item \label{Proposition 2.vi} Stability: If $\Sigma A$ is the suspension of $A \ne \emptyset$, obtained as the quotient space of $A \times [-1,1]$ with $A \times \set{1}$ and $A \times \set{-1}$ collapsed to different points, then $i(\Sigma A) = i(A) + 1$.
\item \label{Proposition 2.vii} Piercing property: If $A$, $C_0$ and $C_1$ are closed, and $\varphi : A \times [0,1] \to C_0 \cup C_1$ is a continuous map such that $\varphi(-u,t) = - \varphi(u,t)$ for all $(u,t) \in A \times [0,1]$, $\varphi(A \times [0,1])$ is closed, $\varphi(A \times \set{0}) \subset C_0$ and $\varphi(A \times \set{1}) \subset C_1$, then $i(\varphi(A \times [0,1]) \cap C_0 \cap C_1) \ge i(A)$.
\item \label{Proposition 2.viii} Neighborhood of zero: If $U$ is a bounded closed symmetric neighborhood of $0$, then $i(\bdry{U}) = \dim W$.
\end{properties}
\end{proposition}

Next we recall the definition and some properties of the pseudo-index of Benci \cite{MR84c:58014}.

\begin{definition}
Let $\A^\ast$ denote the class of symmetric subsets of $W$, let $\M \in \A$ be closed, and let $\Gamma$ denote the group of odd homeomorphisms of $W$ that are the identity outside $E^{-1}(0,b)$. Then the pseudo-index of $A \in \A^\ast$ related to $i$, $\M$, and $\Gamma$ is defined by
\[
i^\ast(A) = \min_{\gamma \in \Gamma}\, i(\gamma(A) \cap \M).
\]
\end{definition}

The following proposition lists some properties of the pseudo-index.

\begin{proposition}[\cite{MR84c:58014}] \label{Proposition 1}
The pseudo-index $i^\ast : \A^\ast \to \N \cup \set{0,\infty}$ has the following properties:
\begin{properties}{i^\ast}
\item \label{Proposition 1.i} If $A \subset B$, then $i^\ast(A) \le i^\ast(B)$.
\item \label{Proposition 1.ii} If $\gamma \in \Gamma$, then $i^\ast(\gamma(A)) = i^\ast(A)$.
\item \label{Proposition 1.iii} If $A$ and $B$ are closed, then $i^\ast(A \cup B) \le i^\ast(A) + i(B)$.
\end{properties}
\end{proposition}

For $j \ge 1$, let
\[
\A^\ast_j = \set{M \in \A^\ast : M \text{ is compact and } i^\ast(M) \ge j}
\]
and set
\begin{equation} \label{5}
c^\ast_j = \inf_{M \in \A^\ast_j}\, \max_{u \in M}\, E(u).
\end{equation}

\begin{theorem} \label{Theorem 1}
Assume \ref{C1} and \ref{C2}. If $0 < c^\ast_{k+1} \le \dotsb \le c^\ast_{k+l} < b$ for some $k \ge 0$ and $l \ge 3$, then $E$ has $(l - 1)/2$ distinct pairs of critical points at levels in $(0,b)$.
\end{theorem}

First we prove a deformation lemma. For $\alpha \le \beta$ in $\R$, let
\[
E_\alpha = \set{u \in W : E(u) \ge \alpha}, \quad E^\beta = \set{u \in W : E(u) \le \beta}, \quad E_\alpha^\beta = \set{u \in W : \alpha \le E(u) \le \beta}.
\]

\begin{lemma} \label{Lemma 1}
Assume \ref{C2} and let $c \in [c^\ast,b)$, $B = K_c \cup M_c \cup -M_c$, and $\delta > 0$. Then there exist $\eps_0 > 0$ and for each $\eps \in (0,\eps_0)$ a map $\eta \in C(W \times [0,1],W)$ satisfying
\begin{enumroman}
\item \label{Lemma 1.i} $\eta(\cdot,0)$ is the identity,
\item \label{Lemma 1.ii} $\eta(\cdot,t)$ is an odd homeomorphism of $W$ for all $t \in [0,1]$,
\item \label{Lemma 1.iii} $\eta(\cdot,t)$ is the identity outside $E_{c - 2 \eps}^{c + 2 \eps} \setminus N_{\delta/3}(B)$ for all $t \in [0,1]$,
\item \label{Lemma 1.iv} $\norm{\eta(u,t) - u} \le \delta/4$ for all $(u,t) \in W \times [0,1]$,
\item \label{Lemma 1.v} $E(\eta(u,\cdot))$ is nonincreasing for all $u \in W$,
\item \label{Lemma 1.vi} $\eta(E^{c + \eps} \setminus N_\delta(B),1) \subset E^{c - \eps}$.
\end{enumroman}
\end{lemma}

\begin{proof}
By \ref{C2}, there exists $\eps_0 > 0$ such that for each $\eps \in (0,\eps_0)$,
\begin{equation} \label{2}
\dnorm{E'(u)} \ge \frac{32 \eps}{\delta} \quad \forall u \in E_{c - 2 \eps}^{c + 2 \eps} \setminus N_{\delta/3}(B).
\end{equation}
Let $V$ be an odd pseudo-gradient vector field for $E$, i.e., a locally Lipschitz continuous mapping from $\set{u \in W : E'(u) \ne 0}$ to $W$ satisfying
\begin{equation} \label{3}
\norm{V(u)} \le \dnorm{E'(u)}, \quad 2 \dualp{E'(u)}{V(u)} \ge \left(\dnorm{E'(u)}\right)^2, \quad V(-u) = -V(u).
\end{equation}
Take an even locally Lipschitz continuous mapping $g : W \to [0,1]$ such that $g = 0$ outside $E_{c - 2 \eps}^{c + 2 \eps} \setminus N_{\delta/3}(B)$ and $g = 1$ on $E_{c - \eps}^{c + \eps} \setminus N_{2 \delta/3}(B)$, and let $\eta(u,t),\, 0 \le t < T(u) \le + \infty$ be the maximal solution of
\[
\dot{\eta} = -4 \eps g(\eta)\, \frac{V(\eta)}{\norm{V(\eta)}^2}, \quad t > 0, \qquad \eta(u,0) = u \in W.
\]
Then
\[
\norm{\eta(u,t) - u} \le 4 \eps \int_0^t \frac{g(\eta(u,\tau))}{\norm{V(\eta(u,\tau))}}\, d\tau \le 8 \eps \int_0^t \frac{g(\eta(u,\tau))}{\dnorm{E'(\eta(u,\tau))}}\, d\tau \le \frac{\delta t}{4}
\]
by \eqref{3} and \eqref{2}, so $\norm{\eta(u,\cdot)}$ is bounded if $T(u) < + \infty$. So $T(u) = + \infty$ and \ref{Lemma 1.i}--\ref{Lemma 1.iv} follow. Since
\begin{equation} \label{4}
\frac{d}{dt}\, (E(\eta(u,t))) = \dualp{E'(\eta)}{\dot{\eta}} = -4 \eps g(\eta)\, \frac{\dualp{E'(\eta)}{V(\eta)}}{\norm{V(\eta)}^2} \le -2 \eps g(\eta) \le 0
\end{equation}
by \eqref{3}, \ref{Lemma 1.v} holds. To see that \ref{Lemma 1.vi} holds, let $u \in E^{c + \eps} \setminus N_\delta(B)$ and suppose that $\eta(u,1) \notin E^{c - \eps}$. Then for all $t \in [0,1]$, $\eta(u,t) \in E_{c - \eps}^{c + \eps}$ by \ref{Lemma 1.v} and $\eta(u,t) \notin N_{2 \delta/3}(B)$ by \ref{Lemma 1.iv}. So $\eta(u,t) \in E_{c - \eps}^{c + \eps} \setminus N_{2 \delta/3}(B)$ and hence $g(\eta(u,t)) = 1$ for all $t \in [0,1]$, so \eqref{4} gives
\[
E(\eta(u,1)) \le E(u) - 2 \eps \le c - \eps,
\]
a contradiction.
\end{proof}

Next we show that if two of the minimax levels defined in \eqref{5} coincide, then there are infinitely many critical points at that level.

\begin{lemma} \label{Lemma 2}
Assume \ref{C1} and \ref{C2}. If $0 < c^\ast_j = c^\ast_{j+1} = c < b$, then $K_c$ is an infinite set.
\end{lemma}

\begin{proof}
If $c \in (0,c^\ast)$, then $E$ satisfies the \PS{c} condition by \ref{C1} and hence the desired conclusion follows from a standard argument. So suppose $c \in [c^\ast,b)$, and let $B = K_c \cup M_c \cup -M_c$. If $K_c$ consists of a finite number of pairs of antipodal points, then \ref{C2} implies that for sufficiently small $\delta > 0$, $N_\delta(B)$ is the disjoint union of $\pm N_\delta(M_c)$ and a finite number of pairs of closed balls centered at antipodal points. So there is an odd continuous map from $N_\delta(B)$ to $\set{\pm 1}$ and hence $i(N_\delta(B)) \le 1$ by \ref{Proposition 2.ii} and \ref{Proposition 2.viii}. We will show that
\begin{equation} \label{8}
i(N_\delta(B)) \ge 2,
\end{equation}
and conclude that $K_c$ is an infinite set.

Let $\eps_0 > 0$, $\eps \in (0,\eps_0)$, and $\eta \in C(W \times [0,1],W)$ be as in Lemma \ref{Lemma 1}. Since $c^\ast_{j+1} = c$, there exists $M \in \A^\ast_{j+1}$ such that $M \subset E^{c + \eps}$ and hence
\begin{equation} \label{6}
j + 1 \le i^\ast(M) \le i^\ast(E^{c + \eps})
\end{equation}
by \ref{Proposition 1.i}. Take $\eps < \min \set{c/2,(b - c)/2}$ and let $\gamma = \eta(\cdot,1)$. Then $\gamma$ is an odd homeomorphism of $W$ by \ref{Lemma 1.ii} and is the identity outside $E^{-1}(0,b)$ by \ref{Lemma 1.iii}, so $\gamma \in \Gamma$. Hence
\begin{equation}
i^\ast(E^{c + \eps}) = i^\ast(\gamma(E^{c + \eps}))
\end{equation}
by \ref{Proposition 1.ii}. By \ref{Lemma 1.vi},
\begin{equation}
\gamma(E^{c + \eps}) = \gamma(E^{c + \eps} \setminus N_\delta(B)) \cup \gamma(N_\delta(B)) \subset E^{c - \eps} \cup \gamma(N_\delta(B)),
\end{equation}
and $\gamma(N_\delta(B))$ is closed since $\gamma$ is a homeomorphism, so
\begin{equation}
i^\ast(\gamma(E^{c + \eps})) \le i^\ast(E^{c - \eps}) + i(\gamma(N_\delta(B)))
\end{equation}
by \ref{Proposition 1.i} and \ref{Proposition 1.iii}. Since $c^\ast_j = c$, $E^{c - \eps} \notin \A^\ast_j$ and hence
\begin{equation}
i^\ast(E^{c - \eps}) \le j - 1,
\end{equation}
and
\begin{equation} \label{7}
i(\gamma(N_\delta(B))) = i(N_\delta(B))
\end{equation}
by \ref{Proposition 2.ii}. Combining \eqref{6}--\eqref{7} gives \eqref{8}.
\end{proof}

We are now ready to prove Theorem \ref{Theorem 1}.

\begin{proof}[Proof of Theorem \ref{Theorem 1}]
We may assume that $0 < c^\ast_{k+1} < \dotsb < c^\ast_{k+l} < b$ in view of Lemma \ref{Lemma 2}. For $j = k + 1,\dots,k + l$, if $E$ satisfies the \PS{c^\ast_j} condition, then $c^\ast_j$ is a critical level of $E$ by a standard argument (see, e.g., Perera et al.\! \cite[Proposition 3.42]{MR2640827}). On the other hand, if $E$ does not satisfy the \PS{c^\ast_j} condition, then $c^\ast_j \in [c^\ast,b)$ by \ref{C1} and $E$ has a \PS{c^\ast_j} sequence with no convergent subsequence. Then $c^\ast_j - c^\ast$ is a critical level of $E$ by \ref{C2}. So $c^\ast_j$ or $c^\ast_j - c^\ast$ is a critical level of $E$ in $(0,b)$ for each $j$ such that $c^\ast_j \ne c^\ast$, and it follows that $E$ has $(l - 1)/2$ distinct critical levels in $(0,b)$.
\end{proof}

Next we prove the following theorem, from which Theorem \ref{Theorem 3} will follow.

\begin{theorem} \label{Theorem 2}
Assume \ref{C1} and \ref{C2}. Let $A_0$ and $B_0$ be symmetric subsets of the unit sphere $S = \set{u \in W : \norm{u} = 1}$ such that $A_0$ is compact, $B_0$ is closed, and
\begin{equation} \label{9}
i(A_0) \ge k + l, \qquad i(S \setminus B_0) \le k
\end{equation}
for some $k \ge 0$ and $l \ge 3$. Assume that there exist $R > r > 0$ such that, setting
\[
A = \set{Ru : u \in A_0}, \qquad B = \set{ru : u \in B_0}, \qquad X = \set{tu : u \in A,\, t \in [0,1]},
\]
we have
\[
\sup_A\, E \le 0 < \inf_B\, E, \qquad \sup_X\, E < b.
\]
Then $E$ has $(l - 1)/2$ distinct pairs of critical points at levels in $(0,b)$.
\end{theorem}

\begin{proof}
We take $\M$ to be the sphere $S_r = \set{u \in W : \norm{u} = r}$, show that
\[
0 < \inf_B\, E \le c^\ast_{k+1} \le \dotsb \le c^\ast_{k+l} \le \sup_X\, E < b,
\]
and apply Theorem \ref{Theorem 1}. We note that $A$ and $S_r \setminus B$ are radially homeomorphic to $A_0$ and $S \setminus B_0$, respectively, and hence
\begin{equation} \label{10}
i(A) \ge k + l, \qquad i(S_r \setminus B) \le k
\end{equation}
by \ref{Proposition 2.ii} and \eqref{9}.

If $M \in \A^\ast_{k+1}$, then \eqref{10} gives
\[
i(S_r \setminus B) \le k < k + 1 \le i^\ast(M) \le i(M \cap S_r)
\]
since the identity is in $\Gamma$, so $M$ intersects $B$ by \ref{Proposition 2.ii}. Hence $c^\ast_{k+1} \ge \inf E(B)$.

For $\gamma \in \Gamma$, consider the continuous map
\[
\varphi : A \times [0,1] \to W, \quad \varphi(u,t) = \gamma(tu).
\]
We have $\varphi(A \times [0,1]) = \gamma(X)$, which is compact. Since $\gamma$ is odd, $\varphi(-u,t) = - \varphi(u,t)$ for all $(u,t) \in A \times [0,1]$ and $\varphi(A \times \set{0}) = \set{\gamma(0)} = \set{0}$. Since $E \le 0$ on $A$, $\restr{\gamma}{A}$ is the identity and hence $\varphi(A \times \set{1}) = A$. Applying the piercing property \ref{Proposition 2.vii} of the index with $C_0 = \set{u \in W : \norm{u} \le r}$ and $C_1 = \set{u \in W : \norm{u} \ge r}$ gives
\[
i(\gamma(X) \cap S_r) = i(\varphi(A \times [0,1]) \cap C_0 \cap C_1) \ge i(A) \ge k + l
\]
by \eqref{10}. Hence $i^\ast(X) \ge k + l$. So $X \in \A^\ast_{k+l}$ and hence $c^\ast_{k+l} \le \sup E(X)$.
\end{proof}

We are now ready to prove Theorems \ref{Theorem 3} and \ref{Theorem 4}.

\begin{proof}[Proof of Theorem \ref{Theorem 3}]
If $C_0 = \emptyset$, then $k = m$ by \eqref{11} and \ref{Proposition 2.i}, and
\[
i(A_0) = i(\varphi(S^N)) \ge i(S^N) = N + 1 = k + N - m + 1
\]
by \ref{Proposition 2.ii} and \ref{Proposition 2.viii}, so the conclusion follows from Theorem \ref{Theorem 2}.

If $C_0 \ne \emptyset$, recall that $\Sigma C_0$ denotes the suspension of $C_0$, which is obtained as the quotient space of $C_0 \times [-1,1]$ with $C_0 \times \set{1}$ and $C_0 \times \set{-1}$ collapsed to different points. Let $\Sigma^{N+1} C_0$ be the $(N + 1)$-fold suspension consisting of points $(v,t_1,\dots,t_{N+1})$, where $v \in C_0$ and $t_j \in [-1,1]$ for $j = 1,\dots,N + 1$, with the appropriate identifications for $t_j = \pm 1$. Set
\[
p_0 = \prod_{l=1}^{N+1} (1 - |t_l|), \qquad p_j = |t_j| \prod_{l=j+1}^{N+1} (1 - |t_l|) \quad \text{for } j = 1,\dots,N, \qquad p_{N+1} = |t_{N+1}|,
\]
let $\set{e_1,\dots,e_{N+1}}$ be the standard unit basis of $\R^{N+1}$, and let
\[
\varpi : \R^{N+1} \setminus \set{0} \to S^N, \quad x \mapsto \frac{x}{|x|}
\]
be the radial projection onto $S^N$. Then
\[
\Sigma^{N+1} C_0 \to A_0, \quad (v,t_1,\dots,t_{N+1}) \mapsto \pi\Bigg(p_0 v + (1 - p_0)\, \varphi\Bigg(\varpi\Bigg(\frac{1}{1 - p_0} \sum_{j=1}^{N+1} p_j \sgn t_j\, e_j\Bigg)\Bigg)\Bigg)
\]
is an odd continuous map, and hence
\[
i(A_0) \ge i(\Sigma^{N+1} C_0) = i(C_0) + N + 1 \ge k + N - m + 1
\]
by \ref{Proposition 2.ii}, \ref{Proposition 2.vi}, and \eqref{11}. So the conclusion follows from Theorem \ref{Theorem 2} again.
\end{proof}

\begin{proof}[Proof of Theorem \ref{Theorem 4}]
We apply Theorem \ref{Theorem 3} with $B_0 = \Psi_{\lambda_{k+1}}$. By \eqref{13}, $i(S \setminus \Psi_{\lambda_{k+1}}) \le k$. For $u \in S$ and $t > 0$,
\[
E(tu) = \frac{t^p}{p} \left(1 - \frac{\lambda}{\Psi(u)}\right) - F(tu).
\]
Since $\Psi(u) > 0$ by $(B_2)$, this gives
\[
\frac{t^p}{p} \left(1 - \frac{\lambda^+}{\Psi(u)}\right) - F(tu) \le E(tu) \le \frac{t^p}{p} \left(1 + \frac{\lambda^-}{\Psi(u)}\right) - F(tu),
\]
where $\lambda^\pm = \max \set{\pm \lambda,0}$. Since $\lambda^+ < \lambda_{k+1}$, the first inequality and $(F_1)$ imply that $\inf E(B) > 0$ if $r > 0$ is sufficiently small. Since $A_0$ is a compact subset of $W \setminus \set{0}$, the second inequality and $(F_2)$ imply that $\sup E(A) \le 0$ if $R > r$ is sufficiently large. By \eqref{18} and \eqref{19}, $\sup E(X) < b$.
\end{proof}

\section{Compactness conditions}

In this section we show that the energy functional
\[
E(u) = \frac{1}{p} \int_\Omega |\nabla u|^p\, dx - \frac{\lambda}{p} \int_\Omega |u|^p\, dx - \frac{1}{p^\ast} \int_\Omega |u|^{p^\ast} dx, \quad u \in W^{1,\,p}_0(\Omega)
\]
satisfies the compactness conditions \ref{C1} and \ref{C2}.

Solutions of the asymptotic equation
\begin{equation} \label{28}
- \Delta_p\, u = |u|^{p^\ast - 2}\, u
\end{equation}
in $\D^{1,\,p}(\R^N)$ coincide with critical points of the functional
\[
E_\infty(u) = \frac{1}{p} \int_{\R^N} |\nabla u|^p\, dx - \frac{1}{p^\ast} \int_{\R^N} |u|^{p^\ast} dx, \quad u \in \D^{1,\,p}(\R^N).
\]
Denote by $\R^N_+ = \set{x = (x_1,\dots,x_N) \in \R^N : x_N > 0}$ the upper-half space in $\R^N$ and by $\D^{1,\,p}_0(\R^N_+)$ the closure of $C^\infty_0(\R^N_+)$ in $\D^{1,\,p}(\R^N)$ after extending by zero on $\R^N \setminus \R^N_+$. Set
\[
c^\ast = \frac{1}{N}\, S_{N,\,p}^{N/p}.
\]

\begin{lemma} \label{nontrivial}
Let $u$ be a nontrivial weak solution of the equation \eqref{28} in $\D^{1,\,p}(\R^N)$ or $\D^{1,\,p}_0(\R^N_+)$. Then
\begin{equation} \label{32}
E_\infty(u) \ge c^\ast.
\end{equation}
If $u \in \D^{1,\,p}_0(\R^N_+)$, then this inequality is strict. If $u$ is sign-changing, then
\begin{equation} \label{33}
E_\infty(u) \ge 2c^\ast.
\end{equation}
\end{lemma}

\begin{proof}
The inequality \eqref{32} follows by testing \eqref{28} with $u$ and using the Sobolev inequality. If $u \in \D^{1,\,p}_0(\R^N_+)$ and equality holds in \eqref{32}, then $u \equiv 0$ by Mercuri and Willem \cite[Theorem 1.1]{MR2644751} (see also Farina et al.\! \cite{MR3989957}). If $u$ is sign-changing, testing \eqref{28} with $u^\pm$ and using the Sobolev inequality gives \eqref{33}.
\end{proof}

Next we prove a global compactness result for problem \eqref{1}.

\begin{theorem} \label{MWglobal}
Let $c \in \R$ and let $\seq{u_n} \subset W^{1,\,p}_0(\Omega)$ be a {\em \PS{c}} sequence for $E$. Then, passing to a subsequence if necessary, there exist a possibly nontrivial solution $u \in W^{1,\,p}_0(\Omega)$ of problem \eqref{1}, $k \in \N \cup \set{0}$, nontrivial solutions $v_i,\, i = 1,\dots,k$ of equation \eqref{28} in $H_i$, where $H_i$ is $\R^N$ or (up to a rotation and a translation) $\R^N_+$, with $v_i \in \D^{1,\,p}(\R^N)$ if $H_i = \R^N$ and $v_i \in \D^{1,\,p}_0(\R^N_+)$ if $H_i = \R^N_+$, and sequences $\seq{y^i_n} \subset \closure{\Omega}$ and $\seq{\eps^i_n} \subset \R_+$ such that
\begin{gather}
\left(\eps^i_n\right)^{-1} \dist{y^i_n}{\bdry{\Omega}} \to \infty \text{ as } n \to \infty \text{ if } H_i = \R^N, \notag\\[10pt]
\left(\eps^i_n\right)^{-1} \dist{y^i_n}{\bdry{\Omega}} \text{ is bounded if } H_i = \R^N_+, \notag\\[10pt]
\label{41} \norm{u_n - u - \sum_{i=1}^k \left(\eps^i_n\right)^{-(N-p)/p} v_i((\cdot - y^i_n)/\eps^i_n)} \to 0 \text{ as } n \to \infty,\\[10pt]
\norm{u_n}^p \to \norm{u}^p + \sum_{i=1}^k \norm{v_i}^p \text{ as } n \to \infty, \notag\\[10pt]
\label{quantum} E(u) + \sum_{i=1}^k E_\infty(v_i) = c.
\end{gather}
\end{theorem}

\begin{proof}
This theorem follows by arguing as in Mercuri and Willem \cite[Theorem 1.2]{MR2644751}. Note that $\seq{u_n}$ is bounded by a standard argument and $k$ is finite by Lemma \ref{nontrivial}. Unlike in \cite{MR2644751}, here $v_i$ may be a solution of \eqref{28} in $\R^N$ or $\R^N_+$ since we make no assumptions on $u_n^-$ (see also Farina et al.\! \cite{MR3989957}).
\end{proof}

Recall that the infimum in \eqref{20} is attained on the functions
\begin{equation} \label{40}
u_{\eps,\,y}(x) = \frac{c_{N,\,p}\, \eps^{(N-p)/p\,(p-1)}}{\big(\eps^{p/(p-1)} + |x - y|^{p/(p-1)}\big)^{(N-p)/p}}, \quad \eps > 0,\, y \in \R^N,
\end{equation}
where the constant $c_{N,\,p} > 0$ is chosen so that
\[
\int_{\R^N} |\nabla u_{\eps,\,y}|^p\, dx = \int_{\R^N} u_{\eps,\,y}^{p^\ast}\, dx = S_{N,\,p}^{N/p}
\]
(see Talenti \cite{MR0463908}). Set
\[
M = \set{u_{\eps,\,y} : \eps > 0,\, y \in \R^N}.
\]
The main result of this section is the following theorem, which implies that $E$ satisfies \ref{C1} and \ref{C2}.

\begin{theorem} \label{PScondition}
The functional $E$ satisfies the {\em \PS{c}} condition for all $c < c^\ast$. If $c^\ast \le c < 2c^\ast$ and $\seq{u_n}$ is a {\em \PS{c}} sequence for $E$ such that $u_n \wto u$ but not strongly, then $u \in K_{c-c^\ast}$ and
\[
\dist{u_n - u}{M} \to 0 \quad \text{or} \quad \dist{u_n - u}{-M} \to 0
\]
for a renamed subsequence. Moreover, setting
\[
M_c = K_{c-c^\ast} + M = \set{u + v : u \in K_{c-c^\ast},\, v \in M},
\]
we have $M_c \subset D \setminus \set{0}$ and $N_\delta(M_c) \cap N_\delta(-M_c) = \emptyset$ for all sufficiently small $\delta > 0$.
\end{theorem}

\begin{proof}
The proof is based on Theorem \ref{MWglobal}. We have
\[
E(u) = E(u) - \frac{1}{p}\, E'(u)\, u = \frac{1}{N} \int_\Omega |u|^{p^\ast} dx
\]
since $E'(u) = 0$, and $E_\infty(v_i) \ge c^\ast$ for $i = 1,\dots,k$ by Lemma \ref{nontrivial}, so \eqref{quantum} yields
\begin{equation} \label{42}
\frac{1}{N} \int_\Omega |u|^{p^\ast} dx + kc^\ast \le c.
\end{equation}
If $c < c^\ast$, this implies $k = 0$, so $u_n \to u$ by \eqref{41}.

Suppose $c^\ast \le c < 2c^\ast$. Then $k \le 1$ by \eqref{42}. If $k = 0$, then $u_n \to u$ as before, so suppose $k = 1$. Then
\[
\norm{u_n - u - \left(\eps^1_n\right)^{-(N-p)/p} v_1((\cdot - y^1_n)/\eps^1_n)} \to 0 \text{ as } n \to \infty
\]
by \eqref{41} and
\begin{equation} \label{43}
E(u) + E_\infty(v_1) = c
\end{equation}
by \eqref{quantum}. If $c = c^\ast$, then $u \equiv 0$ by \eqref{42} and hence $E_\infty(v_1) = c^\ast$ by \eqref{43}. Then either $v_1 \in M$ or $v_1 \in -M$ by Talenti \cite{MR0463908}. If $c^\ast < c < 2c^\ast$, then $E_\infty(v_1) < 2c^\ast$ by \eqref{43} and hence $v_1$ does not change sign by Lemma \ref{nontrivial}. So $v_1$ is a constant sign nontrivial solution of \eqref{28} in $\D^{1,\,p}(\R^N)$ by the Liouville theorem on $\R^N_+$ of Mercuri and Willem \cite[Theorem 1.1]{MR2644751} (see also Farina et al.\! \cite{MR3989957}). Then either $v_1 \in M$ or $v_1 \in -M$ by Sciunzi \cite{MR3459013} and V\'{e}tois \cite{MR3411668} (see also Guedda and V\'{e}ron \cite[Theorem 2.1({\em ii})]{MR964617}). In particular, $E_\infty(v_1) = c^\ast$, so $E(u) = c - c^\ast$ by \eqref{43} and hence $u \in K_{c-c^\ast}$. Since $c - c^\ast < c^\ast$ and hence $K_{c-c^\ast}$ is compact by the first part of the theorem, the rest now follows as in Clapp and Weth \cite[Lemma 9]{MR2122698}.
\end{proof}

\section{Proofs of main theorems} \label{Proofs}

In this section we prove Theorems \ref{Theorem 5} and \ref{Theorem 6} by applying Theorem \ref{Theorem 4} with $D = \D^{1,\,p}(\R^N)$, $W = W^{1,\,p}_0(\Omega)$, and the operators $\Ap, \Bp, f \in C(W^{1,\,p}_0(\Omega),W^{-1,\,p'}(\Omega))$ given by
\[
\dualp{\Ap[u]}{v} = \int_\Omega |\nabla u|^{p-2}\, \nabla u \cdot \nabla v\, dx, \hquad \dualp{\Bp[u]}{v} = \int_\Omega |u|^{p-2}\, uv\, dx, \hquad \dualp{f(u)}{v} = \int_\Omega |u|^{p^\ast - 2}\, uv\, dx
\]
for $u, v \in W^{1,\,p}_0(\Omega)$. It is easily seen that $(A_1)$--$(F_2)$ hold, and $E$ satisfies \ref{C1} and \ref{C2} with
\[
c^\ast = \frac{1}{N}\, S_{N,\,p}^{N/p}, \qquad b = \frac{2}{N}\, S_{N,\,p}^{N/p}
\]
by Theorem \ref{PScondition}.

To construct the compact symmetric set $C_0$, we assume without loss of generality that $0 \in \Omega$. Fix $0 < \delta_0 < \dist{0}{\bdry{\Omega}}$ and let $\pi : W^{1,\,p}_0(\Omega) \setminus \set{0} \to S,\, u \mapsto u/\norm{u}$ be the radial projection onto $S$.

\begin{lemma} \label{Lemma 3}
For each $k \ge 1$ and sufficiently small $0 < \delta < \delta_0$, there exists a compact symmetric subset $C_{k,\delta}$ of $S$ with $i(C_{k,\delta}) \ge k$ such that $v = 0$ on $B_{3 \delta/4}(0)$ for all $v \in C_{k,\delta}$ and
\begin{equation} \label{34}
\sup_{v \in C_{k,\delta},\, s \ge 0}\, E(sv) \le \begin{cases}
0 & \text{if } \lambda > \lambda_k\\[7.5pt]
a \delta^{(N-p)\,N/p} & \text{if } \lambda = \lambda_k
\end{cases}
\end{equation}
for some constant $a > 0$.
\end{lemma}

\begin{proof}
We have $\lambda_k = \cdots = \lambda_l < \lambda_{l+1}$ for some $l \ge k$. By Degiovanni and Lancelotti \cite[Theorem 2.3]{MR2514055}, the sublevel set $\Psi^{\lambda_l}$ has a compact symmetric subset $C_l$ of index $l$ that is bounded in $L^\infty(\Omega) \cap C^{1,\,\alpha}_\loc(\Omega)$. Let $\xi : [0,\infty) \to [0,1]$ be a smooth function such that $\xi(s) = 0$ for $s \le 3/4$ and $\xi(s) = 1$ for $s \ge 1$. For $u \in C_l$, set
\[
u_\delta(x) = \xi\bigg(\frac{|x|}{\delta}\bigg)\, u(x), \qquad v = \pi(u_\delta),
\]
and let
\[
C_{k,\delta} = \set{v : u \in C_l}.
\]
Since $C_l$ is a compact symmetric set and $u \mapsto v$ is an odd continuous map of $C_l$ onto $C_{k,\delta}$, $C_{k,\delta}$ is also a compact symmetric set and
\[
i(C_{k,\delta}) \ge i(C_l) = l \ge k
\]
by \ref{Proposition 2.ii}.

Let $u \in C_l$, $v = \pi(u_\delta)$, and $s \ge 0$. Since $u$ is bounded in $C^1(B_{\delta_0}(0))$ and belongs to $\Psi^{\lambda_l} = \Psi^{\lambda_k}$,
\[
\int_\Omega |\nabla u_\delta|^p\, dx \le \int_{\Omega \setminus B_\delta(0)} |\nabla u|^p\, dx + \int_{B_\delta(0)} \left(|\nabla u|^p + a_1\, \delta^{-p}\, |u|^p\right) dx \le 1 + a_2\, \delta^{N-p}
\]
and
\[
\int_\Omega |u_\delta|^p\, dx \ge \int_{\Omega \setminus B_\delta(0)} |u|^p\, dx = \int_\Omega |u|^p\, dx - \int_{B_\delta(0)} |u|^p\, dx \ge \frac{1}{\lambda_k} - a_3\, \delta^N
\]
for some constants $a_1, a_2, a_3 > 0$. So
\begin{equation} \label{30}
\int_\Omega |v|^p\, dx = \frac{\dint_\Omega |u_\delta|^p\, dx}{\dint_\Omega |\nabla u_\delta|^p\, dx} \ge \frac{1}{\lambda_k} - a_4\, \delta^{N-p}
\end{equation}
for some constant $a_4 > 0$ and hence
\begin{equation} \label{31}
\int_\Omega |v|^{p^\ast} dx \ge a_5
\end{equation}
for some constant $a_5 > 0$. We have
\[
E(sv) = \frac{s^p}{p} \int_\Omega |\nabla v|^p\, dx - \frac{\lambda s^p}{p} \int_\Omega |v|^p\, dx - \frac{s^{p^\ast}}{p^\ast} \int_\Omega |v|^{p^\ast} dx,
\]
and maximizing the right-hand side over $s \ge 0$ gives
\[
E(sv) \le \frac{1}{N}\: Q(v)^{N/p},
\]
where
\[
Q(v) = \frac{\left(\dint_\Omega |\nabla v|^p\, dx - \lambda \dint_\Omega |v|^p\, dx\right)^+}{\left(\dint_\Omega |v|^{p^\ast} dx\right)^{p/p^\ast}}.
\]
By \eqref{30} and \eqref{31},
\[
Q(v) \le a_6 \left(1 - \frac{\lambda}{\lambda_k} + a_7\, \delta^{N-p}\right)^+
\]
for some constants $a_6, a_7 > 0$, so \eqref{34} follows for sufficiently small $\delta > 0$.
\end{proof}

Next we construct the odd continuous map $\varphi$.

\begin{lemma} \label{Lemma 4}
For each $\eps > 0$ and $0 < \delta < \delta_0$ with $\eps \ll \delta$, there exists an odd continuous map $\varphi_{\eps,\delta} : S^N \to S$ such that $w = 0$ on $\Omega \setminus B_{3 \delta/4}(0)$ for all $w \in \varphi_{\eps,\delta}(S^N)$ and
\begin{equation} \label{26}
\sup_{w \in \varphi_{\eps,\delta}(S^N),\, t \ge 0}\, E(tw) \le \begin{cases}
\ds{\frac{2}{N}\, S_{N,\,p}^{N/p} \left[1 + a_1 \bigg(\frac{\eps}{\delta}\bigg)^{(N-p)/(p-1)} - a_2\, \eps^p\right]^{N/p}} & \text{if } N > p^2\\[20pt]
\ds{\frac{2}{p^2}\, S_{N,\,p}^p\, \Bigg[1 + a_1 \bigg(\frac{\eps}{\delta}\bigg)^p - a_2\, \eps^p \log \bigg(\frac{\delta}{\eps}\bigg)\Bigg]^p} & \text{if } N = p^2
\end{cases}
\end{equation}
for some constants $a_1, a_2 > 0$.
\end{lemma}

\begin{proof}
Referring to \eqref{40}, let
\[
u_\eps(x) = u_{\eps,\,0}(x) = \frac{c_{N,\,p}\, \eps^{(N-p)/p\,(p-1)}}{\big(\eps^{p/(p-1)} + |x|^{p/(p-1)}\big)^{(N-p)/p}}, \quad \eps > 0.
\]
Let $\zeta : [0,\infty) \to [0,1]$ be a smooth function such that $\zeta(s) = 1$ for $s \le 1/8$ and $\zeta(s) = 0$ for $s \ge 1/4$, and set
\[
u_{\eps,\delta}(x) = \zeta\bigg(\frac{|x|}{\delta}\bigg)\, u_\eps(x), \quad \eps > 0,\, 0 < \delta < \delta_0.
\]
We have the estimates
\begin{gather}
\label{21} \int_\Omega |\nabla u_{\eps,\delta}|^p\, dx \le S_{N,\,p}^{N/p} \left[1 + a_3 \bigg(\frac{\eps}{\delta}\bigg)^{(N-p)/(p-1)}\right],\\[15pt]
\left(\int_\Omega u_{\eps,\delta}^{p^\ast}\, dx\right)^{p/p^\ast} \ge S_{N,\,p}^{N/p^\ast} \left[1 - a_4 \bigg(\frac{\eps}{\delta}\bigg)^{N/(p-1)}\right],\\[15pt]
\label{22} \int_\Omega u_{\eps,\delta}^p\, dx \ge \begin{cases}
\ds{\eps^p \left[a_5 - a_6 \bigg(\frac{\eps}{\delta}\bigg)^{(N-p^2)/(p-1)}\right]} & \text{if } N > p^2\\[20pt]
\ds{\eps^p\, \Bigg[a_5 \log \bigg(\frac{\delta}{\eps}\bigg) - a_6\Bigg]} & \text{if } N = p^2
\end{cases}
\end{gather}
for some constants $a_3, a_4, a_5, a_6 > 0$ (see, e.g., Degiovanni and Lancelotti \cite[Lemma 3.1]{MR2514055}).

Let $S^{N-1}$ be the unit sphere in $\R^N$, let
\[
S^N_+ = \big\{x = (x' \sqrt{1 - s^2},s) : x' \in S^{N-1},\, s \in [0,1]\big\}
\]
be the upper hemisphere in $\R^{N+1}$, and define a continuous map $\varphi_{\eps,\delta} : S^N_+ \to S$ by
\[
\varphi_{\eps,\delta}(x) = \pi\big(u_{\eps,\delta}(\cdot - (1 - s)\, x'/2) - (1 - s)\, u_{\eps,\delta}(\cdot + x'/2)\big).
\]
Let $w = \varphi_{\eps,\delta}(x)$ and note that $w = 0$ outside $B_{3 \delta/4}(0)$. For $t \ge 0$,
\[
E(tw) = \frac{t^p}{p} \int_\Omega |\nabla w|^p\, dx - \frac{\lambda t^p}{p} \int_\Omega |w|^p\, dx - \frac{t^{p^\ast}}{p^\ast} \int_\Omega |w|^{p^\ast} dx,
\]
and maximizing the right-hand side over $t \ge 0$ gives
\[
E(tw) \le \frac{1}{N}\: Q(w)^{N/p},
\]
where
\[
Q(w) = \frac{\dint_\Omega |\nabla w|^p\, dx - \lambda \dint_\Omega |w|^p\, dx}{\left(\dint_\Omega |w|^{p^\ast} dx\right)^{p/p^\ast}}.
\]
Noting that $Q(w) = Q(u_{\eps,\delta}(\cdot - (1 - s)\, x'/2) - (1 - s)\, u_{\eps,\delta}(\cdot + x'/2))$ and $u_{\eps,\delta}(\cdot - (1 - s)\, x'/2)$ and $u_{\eps,\delta}(\cdot + x'/2)$ have disjoint supports gives
\[
Q(w) = \frac{1 + (1 - s)^p}{[1 + (1 - s)^{p^\ast}]^{p/p^\ast}}\; Q(u_{\eps,\delta}),
\]
and maximizing the right-hand side over $s \in [0,1]$ gives
\[
Q(w) \le 2^{p/N} Q(u_{\eps,\delta}).
\]
So
\begin{equation} \label{25}
E(tw) \le \frac{2}{N}\: Q(u_{\eps,\delta})^{N/p}.
\end{equation}
Since $\varphi_{\eps,\delta}$ is odd on $S^{N-1}$ and $E$ is an even functional, $\varphi_{\eps,\delta}$ can be extended to an odd continuous map from $S^N$ to $S$ such that \eqref{25} holds for all $w \in \varphi_{\eps,\delta}(S^N)$ and $t \ge 0$. Since $\eps \ll \delta$,
\[
Q(u_{\eps,\delta}) \le \begin{cases}
\ds{S + a_3 \bigg(\frac{\eps}{\delta}\bigg)^{(N-p)/(p-1)} - a_4\, \eps^p} & \text{if } N > p^2\\[20pt]
\ds{S + a_3 \bigg(\frac{\eps}{\delta}\bigg)^p - a_4\, \eps^p \log \bigg(\frac{\delta}{\eps}\bigg)} & \text{if } N = p^2
\end{cases}
\]
for some constants $a_3, a_4 > 0$ by \eqref{21}--\eqref{22}, so \eqref{26} follows.
\end{proof}

We are now ready to prove Theorems \ref{Theorem 5} and \ref{Theorem 6}.

\begin{proof}[Proof of Theorem \ref{Theorem 5}]
Let $0 < \delta < \delta_0$. Lemma \ref{Lemma 4} gives an odd continuous map $\varphi_{\eps,\delta} : S^N \to S$ satisfying
\begin{equation} \label{27}
\sup_{w \in \varphi_{\eps,\delta}(S^N),\, t \ge 0}\, E(tw) < \frac{2}{N}\, S_{N,\,p}^{N/p}
\end{equation}
for all sufficiently small $\eps > 0$.

\ref{Theorem 5.i} If $0 < \lambda < \lambda_1$, Theorem \ref{Theorem 4} with $k = m = 0$ and $C_0 = \emptyset$ gives $N/2$ distinct pairs of nontrivial solutions satisfying \eqref{23}. If $\lambda_k < \lambda < \lambda_{k+1}$, Lemma \ref{Lemma 3} gives a compact symmetric subset $C_{k,\delta}$ of $S$ with $i(C_{k,\delta}) \ge k$ satisfying
\begin{equation}
\sup_{v \in C_{k,\delta},\, s \ge 0}\, E(sv) = 0
\end{equation}
when $\delta$ is sufficiently small. Let $v \in C_{k,\delta}$, $w \in \varphi_{\eps,\delta}(S^N)$, and $s, t \ge 0$. Since $v = 0$ on $B_{3 \delta/4}(0)$ and $w = 0$ outside $B_{3 \delta/4}(0)$, $\varphi_{\eps,\delta}(S^N) \subset S \setminus C_{k,\delta}$ and
\begin{equation} \label{29}
E(sv + tw) = E(sv) + E(tw).
\end{equation}
It follows from \eqref{27}--\eqref{29} that
\[
\sup_{v \in C_{k,\delta},\, w \in \varphi_{\eps,\delta}(S^N),\, s, t \ge 0}\, E(sv + tw) < \frac{2}{N}\, S_{N,\,p}^{N/p}.
\]
So Theorem \ref{Theorem 4} with $m = 0$ again gives $N/2$ distinct pairs of nontrivial solutions satisfying \eqref{23}.

\ref{Theorem 5.ii} If $\lambda = \lambda_1$, then $\lambda < \lambda_2$ since the first eigenvalue is simple, so Theorem \ref{Theorem 4} with $k = m = 1$ and $C_0 = \emptyset$ gives $(N - 1)/2$ distinct pairs of nontrivial solutions satisfying \eqref{23}.

\ref{Theorem 5.iii} If $\lambda_{k-m} < \lambda = \lambda_{k-m+1} = \cdots = \lambda_k < \lambda_{k+1}$, where $k > m \ge 1$, Lemma \ref{Lemma 3} gives a compact symmetric subset $C_{k-m,\delta}$ of $S$ with $i(C_{k-m,\delta}) \ge k - m$ satisfying
\[
\sup_{v \in C_{k-m,\delta},\, s \ge 0}\, E(sv) = 0
\]
when $\delta$ is sufficiently small. As in the proof of part \ref{Theorem 5.i}, $\varphi_{\eps,\delta}(S^N) \subset S \setminus C_{k-m,\delta}$ and
\[
\sup_{v \in C_{k-m,\delta},\, w \in \varphi_{\eps,\delta}(S^N),\, s, t \ge 0}\, E(sv + tw) < \frac{2}{N}\, S_{N,\,p}^{N/p}.
\]
So Theorem \ref{Theorem 4} gives $(N - m)/2$ distinct pairs of nontrivial solutions satisfying \eqref{23}.
\end{proof}

\begin{proof}[Proof of Theorem \ref{Theorem 6}]
The case where $0 < \lambda < \lambda_1$ or $\lambda_k < \lambda < \lambda_{k+1}$ for some $k \ge 1$ is covered in Theorem \ref{Theorem 5} \ref{Theorem 5.i}, so we assume that $\lambda = \lambda_k < \lambda_{k+1}$ for some $k \ge 1$. Lemma \ref{Lemma 3} gives a compact symmetric subset $C_{k,\delta}$ of $S$ with $i(C_{k,\delta}) \ge k$ satisfying
\begin{equation} \label{35}
\sup_{v \in C_{k,\delta},\, s \ge 0}\, E(sv) \le a \delta^{(N-p)\,N/p}
\end{equation}
for some constant $a > 0$ when $0 < \delta < \delta_0$ is sufficiently small. Lemma \ref{Lemma 4} gives an odd continuous map $\varphi_{\eps,\delta} : S^N \to S$ satisfying
\begin{equation}
\sup_{w \in \varphi_{\eps,\delta}(S^N),\, t \ge 0}\, E(tw) \le \frac{2}{N}\, S_{N,\,p}^{N/p} \left[1 + a_1 \bigg(\frac{\eps}{\delta}\bigg)^{(N-p)/(p-1)} - a_2\, \eps^p\right]^{N/p}
\end{equation}
for some constants $a_1, a_2 > 0$ for all sufficiently small $\eps > 0$. Let $v \in C_{k,\delta}$, $w \in \varphi_{\eps,\delta}(S^N)$, and $s, t \ge 0$. Since $v = 0$ on $B_{3 \delta/4}(0)$ and $w = 0$ outside $B_{3 \delta/4}(0)$, $\varphi_{\eps,\delta}(S^N) \subset S \setminus C_{k,\delta}$ and
\begin{equation}
E(sv + tw) = E(sv) + E(tw).
\end{equation}
Take $\delta = \eps^\alpha$ with $0 < \alpha < 1$. Then
\begin{equation} \label{36}
\delta^{(N-p)\,N/p} = \eps^{p+[(N-p)\,N\alpha-p^2]/p}, \qquad \bigg(\frac{\eps}{\delta}\bigg)^{(N-p)/(p-1)} = \eps^{p+[N-p^2-(N-p)\,\alpha]/(p-1)}.
\end{equation}
Since $N^2/(N + 1) > p^2$,
\[
0 < \frac{p^2}{(N - p)\, N} < \frac{N - p^2}{N - p} < 1,
\]
so we can take $p^2/(N - p)\, N < \alpha < (N - p^2)/(N - p)$, combine \eqref{35}--\eqref{36}, and take $\eps$ sufficiently small to get
\[
\sup_{v \in C_{k,\delta},\, w \in \varphi_{\eps,\delta}(S^N),\, s, t \ge 0}\, E(sv + tw) < \frac{2}{N}\, S_{N,\,p}^{N/p}.
\]
So Theorem \ref{Theorem 4} with $m = 0$ gives $N/2$ distinct pairs of nontrivial solutions satisfying \eqref{24}.
\end{proof}

\def\cdprime{$''$}

\end{document}